\documentclass[12pt]{amsart}
\usepackage{amscd,amssymb}
\usepackage[graph,frame,poly,arc]{xy}  
\usepackage[plainpages,backref,urlcolor=blue]{hyperref}

\theoremstyle{plain}
\newtheorem{thm}[subsection]{Theorem}
\newtheorem{lem}[subsection]{Lemma}
\newtheorem{prop}[subsection]{Proposition}
\newtheorem{cor}[subsection]{Corollary}

\theoremstyle{definition}
\newtheorem{rk}[subsection]{Remark}
\newtheorem{definition}[subsection]{Definition}
\newtheorem{ex}[subsection]{Example}

\numberwithin{equation}{section}
\newcommand{\A}{{\mathcal A}}

\newcommand{\C}{\mathbb{C}}
\newcommand{\PP}{\mathbb{P}}

\DeclareMathOperator{\rank}{rank}

\begin{document}

\title [Freeness  versus maximal global Tjurina number for plane curves]
{Freeness  versus maximal global Tjurina number for plane curves }

\author[Alexandru Dimca]{Alexandru Dimca$^1$}
\address{Univ. Nice Sophia Antipolis, CNRS,  LJAD, UMR 7351, 06100 Nice, France. }
\email{dimca@unice.fr}

\thanks{$^1$ Partially supported by Institut Universitaire de France.}



\begin{abstract} 
We give a characterization of nearly free plane curves in terms of their global Tjurina numbers, 
similar to the characterization of free curves as curves with a maximal Tjurina number, given by A. A. du Plessis and C.T.C. Wall.
It is also shown that an irreducible plane curve having a 1-dimensional symmetry is nearly free. 
A new numerical characterization of free curves and a simple characterization of nearly free curves in terms of their syzygies conclude this note.

\end{abstract}
 
\maketitle


\section{Introduction} 
This note is inspired and should be regarded as a modest continuation of the beautiful paper \cite{duPCTC} by A. A. du Plessis and C.T.C. Wall. We start by quoting a part of the main result of this paper. 
Let $S=\C[x,y,z]$ be the graded polynomial ring in three variables $x,y,z$ and let $C:f=0$ be a reduced curve of degree $d$ in the complex projective plane $\PP^2$. The minimal degree of a Jacobian relation for $f$ is the integer $mdr(f)$
defined to be the smallest integer $m\geq 0$ such that there is a nontrivial relation
\begin{equation}
\label{rel_m}
 af_x+bf_y+cf_z=0
\end{equation}
among the partial derivatives $f_x, f_y$ and $f_z$ of $f$ with coefficients $a,b,c$ in $S_m$, the vector space of  homogeneous polynomials of degree $m$. When $mdr(f)=0$, then $C$ is a union of lines passing through one point, a situation easy to analyse. We assume from now on that 
$$ mdr(f)\geq 1.$$
Denote by $\tau(C)$ the global Tjurina number of the curve $C$, which is the sum of the Tjurina numbers of the singular points of $C$.
Then the result of du Plessis and Wall referred to above is the following, see Theorem 3.2 in \cite{duPCTC}.
\begin{thm}
\label{thmCTC}
For a positive integer $r$, define two integers by 
$$\tau(r)_{min}=(d-1)(d-r-1) \text { and } \tau(r)_{max}=  (d-1)(d-r-1)+r^2.$$
If $r=mdr(f) <d/2$, then one has
$\tau(r)_{min} \leq \tau(C) \leq \tau(r)_{max}.$
Moreover, if $d$ is even and $r=d/2$, then
$ \tau(r)_{min} \leq \tau(C) \leq \tau(r)_{max}-1.$
\end{thm}

At the end of the proof of this result, the authors state the following very interesting consequence, in a rather hidden way (which prevented us for noticing it for some time).
\begin{cor}
\label{corCTC}
If $r=mdr(f) <d/2$, then one has
$$ \tau(C) =\tau(r)_{max}$$
if and only if $C:f=0$ is a free curve.
\end{cor}
The basic properties of the free curves are reviewed in the next section, for now we just say that this is the same as asking the surface singularity given by the cone over $C$ to be a free divisor germ in $(\C^3,0)$ in the sense of K. Saito, who introduced the important notion of free divisor in \cite{KS}.

Recently a related notion, namely that of a nearly free curve, was introduced by G. Sticlaru and the author in  \cite{DStNF}, motivated by the study of rational cuspidal curves. The main result of this note is the following.

\begin{thm}
\label{thm1}
If $r=mdr(f)  \leq d/2$, then one has
$$ \tau(C) = \tau(r)_{max}-1$$
if and only if $C$ is a nearly free curve.
\end{thm}

The case $r=1$ deserves special attention, see also \cite{duPCTC2}, Prop. 1.1 and Prop. 1.3.

\begin{cor}
\label{cor1}
\noindent (i) One has $mdr(f)=1$ if and only if $C$ admits a 1-dimensional symmetry, i.e. $C$ admits a
1-dimensional algebraic subgroup of $PGL_2(\C)$ as automorphism group.

\noindent (ii) If $mdr(f)=1$, then $C$ is either free or nearly free.
If in addition $C$ is irreducible, then $C$ is nearly free.
\end{cor}

To state our final result, we recall  some definitions, see \cite{DStFD}. We denote by $J_f$ the Jacobian ideal of $f$, i.e. the homogeneous ideal in $S$ spanned by $f_x,f_y,f_z$, and  by $M(f)=S/J_f$ the corresponding graded ring, called the Jacobian (or Milnor) algebra of $f$.

\begin{definition}
\label{def}

\noindent (i) the {\it coincidence threshold} 
$$ct(f)=\max \{q:\dim M(f)_k=\dim M(f_s)_k \text{ for all } k \leq q\},$$
with $f_s$  a homogeneous polynomial in $S$ of the same degree $d$ as $f$ and such that $C_s:f_s=0$ is a smooth curve in $\PP^2$.

\noindent (ii) the {\it stability threshold} 
$st(f)=\min \{q~~:~~\dim M(f)_k=\tau(C) \text{ for all } k \geq q\}.$

\end{definition}

It is clear that one has
\begin{equation} 
\label{REL}
ct(f) \geq mdr(f)+d-2,
\end{equation} 
with equality for $mdr(f) <d-1.$
It is interesting that the freeness of the plane curve $C$ can be characterized in terms of these invariants. The first part of the following result was proved in \cite{DStFD}, while the second part
was conjectured in \cite{DStFD} in a weaker form and it is proved below using Theorem \ref{thm1} and additional results from du Plessis and Wall paper \cite{duPCTC}.

\begin{thm}
\label{thmCONJ} 
 (i) For a reduced free plane curve $C:f=0$ of degree $d$ one has
$$ct(f)+st(f)=T,$$
where $T=3(d-2).$

\medskip

\noindent (ii) Conversely, suppose that the reduced  plane curve $C:f=0$ of degree $d$ satisfies
$$ct(f)+st(f) \leq T+1.$$
Then $C$ is free.
\end{thm}

\begin{cor}
\label{cor21}
\noindent (i) For a reduced curve $C:f=0$ one has $ct(f)+st(f) \geq T$ and the equality holds if and only if $C$ is free.

\noindent (ii) For a reduced non free curve $C:f=0$ one has $ct(f)+st(f) \geq T+2$ and the equality holds if and only if $C$ is nearly free.
\end{cor}

In the second section we collect some basic facts on free and nearly free curves. Then we give a new proof to Corollary \ref{corCTC}, perhaps in a language more familiar to people in algebraic geometry and commutative algebra than the proofs in  \cite{duPCTC}. We also state a result related to Terao's conjecture in the case of line arrangements, see Corollary \ref{cor2}.

Then, in the third section, we use exactly the same approach as in the proof of Corollary \ref{corCTC}, in addition to a key exact sequence \eqref{es} introduced by du Plessis and Wall,  to prove the new results, namely Theorem \ref{thm1} and the part (ii) in Theorem \ref{thmCONJ} and Corollary \ref{cor1}.
As a byproduct, we obtain in the last section a new, very simple characterization of nearly free curves, see Theorem \ref{thm2} (ii).

\section{Free and nearly free plane curves}

 Let $I_f$ denote the saturation of the ideal $J_f$ with respect to the maximal ideal $(x,y,z)$ in $S$ and let $N(f)=I_f/J_f$ be the corresponding quotient.

 Consider the graded $S-$submodule $AR(f) \subset S^{3}$ of {\it all relations} involving the derivatives of $f$, namely
$$\rho=(a,b,c) \in AR(f)_m$$
if and only if  $af_x+bf_y+cf_z=0$ and $a,b,c$ are in $S_m$. We set $ar(f)_k=\dim AR(f)_k$,   $m(f)_k=\dim M(f)_k$ and $n(f)_k=\dim N(f)_k$ for any integer $k$.

\begin{definition} The curve $C:f=0$ is a { free divisor} if the following  equivalent conditions hold.

\begin{enumerate}

\item $N(f)=0$, i.e. the Jacobian ideal is saturated.

\item The minimal resolution of the Milnor algebra $M(f)$ has the following  form
$$0 \to S(-d_1-d+1) \oplus S(-d_2-d+1) \to S^3(-d+1) \xrightarrow{(f_x,f_y,f_z)}  S$$
for some positive integers $d_1, d_2$.
\item The graded $S$-module $AR(f)$ is free of rank 2, i.e. there is an isomorphism 
$$AR(f)=S(-d_1) \oplus S(-d_2)$$
for some positive integers $d_1, d_2$.
\end{enumerate}

\end{definition}
When $C$ is a free divisor, the integers $d_1 \leq d_2$ are called the {exponents} of $C$.  They satisfy the relations 
\begin{equation}
\label{free1}
 d_1+d_2=d-1 \text{ and } \tau(C)=(d-1)^2 - d_1d_2,
\end{equation}
where $\tau(C)$ is the total Tjurina number of $C$, see for instance \cite{DS14}, \cite{DStFD}.

Consider the rank two vector bundle $T\langle C\rangle=Der(-logC)$ of logarithmic vector fields  along $C$, which is the coherent sheaf associated to the graded $S$-module $AR(f)(1)$.
Using the results in the third section of \cite{DS14}, for any integer $k$ one has
\begin{equation}
\label{chi}
\chi(T\langle C\rangle(k))=3{k+3 \choose 2}-{d+k+2 \choose 2} +\tau(C).
\end{equation}
Moreover, one has the following for $E=T\langle C\rangle$ and any integer $k$, see \cite{DS14}, \cite{Se}.
\begin{equation}
\label{chi2}
 h^0((E(k))=ar(f)_{k+1}, \  \  h^1((E(k))=n(f)_{d+k} \text{ and
}  h^2((E(k))=ar(f)_{d-5-k}.
\end{equation}
Note that $C$ is free if and only if the vector bundle $T\langle C\rangle$ splits as a direct sum of two line bundles on $\PP^2$.
The definition of a nearly free curve is more subtle, see \cite{DStNF}, combined with Remark 5.2 and Theorem 5.3 (in fact the version of it corresponding to curves) in \cite{DStFS}.
\begin{definition} 
\label{def2}

The curve $C:f=0$ is a { nearly free divisor} 
if the following  equivalent conditions hold.

\begin{enumerate}

\item $N(f) \ne 0$ and $n(f)_k \leq 1$ for any $k$.

\item The Milnor algebra $M(f)$ has a minimal resolution of the form
$$0 \to S(-d-d_2) \to S(-d-d_1+1) \oplus S^2(-d-d_2+1) \to S^3(-d+1) \xrightarrow{(f_0,f_1,f_2)}  S$$
for some integers $1 \leq d_1 \leq d_2$, called the exponents of $C$.

\item There are 3 syzygies $\rho_1$, $\rho_2$, $\rho_3$ of degrees $d_1$, $d_2=d_3=d-d_1$ which form a minimal system of generators for the first syzygies module $AR(f)$.

\end{enumerate}

\end{definition}

If $C:f=0$ is nearly free, then the exponents $d_1 \leq d_2$ satisfy 
\begin{equation}
\label{nfree1}
 d_1+d_2=d \text{ and } ,\tau(C)=(d-1)^2-d_1(d_2-1)-1,
\end{equation}
see \cite{DStNF}.
For both a free and a nearly free curve $C:f=0$, it is clear that $mdr(f)=d_1$.

\subsection{A new proof for Corollary \ref{corCTC}} \label{proofCTC}

If $C$ is free and $d_1=mdr(f)=r$, it follows from \eqref{free1} that $d_2=d-r-1$ and
$$\tau(C)=(d-1)^2-r(d-r-1)=(d-1)(d-r-1)+r^2.$$
The converse implication is more involved. To estimate the dimension $ar(f)_{d-r-1}$, we use the formula \eqref{chi} given above for $\chi(T\langle C\rangle(k))$, with $k=d-r-2$. Since
$$h^2((T\langle C\rangle(k))=ar(f)_{d-5-k}=ar(f)_{r-3}=0,$$
as $r$ is the minimal degree of an element in $AR(f)$, it follows that
$$ar(f)_{d-r-1}-n(f)_{2d-r-2}=3{d-r+1 \choose 2}-{2d-r \choose 2}+\tau(C).$$
A direct computation using $\tau(C)=(d-1)(d-r-1)+r^2$ yields
$$ar(f)_{d-r-1}-n(f)_{2d-r-2}={d-2r+1 \choose 2}+1.$$
Note that ${d-2r+1 \choose 2}=\dim S_{d-2r-1}$. If we denote by $\rho_1 \in AR(f)$ the relation of minimal degree $r$, the vector space $S_{d-2r-1}\rho_1$ is contained in $AR(f)_{d-r-1}$ and has the dimension ${d-2r+1 \choose 2}$. It follows that there is at least one relation 
$$\rho_2 \in AR(f)_{d-r-1} \setminus S_{d-2r-1}\rho_1.$$
 Then Lemma 1.1 in \cite{ST} implies that $C$ is a free divisor with exponents $d_1=r$ and $d_2=d-1-r$.

\subsection{An application to Terao's Conjecture} \label{Terao}

H. Terao has conjectured that if $\A$ and $\A'$ are hyperplane arrangements in $\PP^n$  with isomorphic intersection lattices $L(\A) =L(\A')$, and if $\A$ is free, then $\A'$ is also free, see  for details
\cite{OT}, \cite{Yo} as well as \cite{HS} for the case $n=2$.
Using Theorem \ref{thmCTC} and Corollary \ref{corCTC}, we get the following partial positive answer in the case of line arrangements.

\begin{cor}
\label{cor2} Let $\A:f=0$ and $\A':f'=0$ be two line arrangements in $\PP^2$  with isomorphic intersection lattices $L(\A) =L(\A')$. Assume $\A$ consists of $d\geq 3$ lines and consider its total Tjurina number
$$\tau(C)=\sum_p(n(p)-1)^2,$$
where the sum is over all multiple points $p$ of $\A$ and $n(p)$ denotes the multiplicity of $\A$ at $p$. If $\A$ is free, then there is  a unique integer $r\geq 0$ such that $r<d/2$ and
$$\tau(\A) =(d-1)(d-r-1)+r^2.$$
If this integer $r$ satisfies $r\leq \sqrt {d-3}$, then the line arrangement $\A'$ is also free.

\end{cor}

\proof A  computation shows that for $s \leq \sqrt {d-3}$, the intervals
$[\tau(s)_{min},\tau(s)_{max}]$ and $[\tau(s-1)_{min},\tau(s-1)_{max}]$ are disjoint, since
$$\tau(s)_{max} <\tau(s-1)_{min}.$$
It follows that each value $\tau(s)_{max}$ uniquely determines the corresponding $s$ when $s\leq \sqrt {d-3}$. 
Moreover, $\tau(\A)=\tau(\A')$, since this number depends only on the intersection lattice  $L(\A) =L(\A')$. It follows that $mdr(f')=mdr(f)=r$ and hence $\A'$ is free by Corollary \ref{corCTC}.

\endproof

\begin{ex}
\label{ex12} When $d=12$, the intervals $[\tau(r)_{min},\tau(r)_{max}]$ are listed in \cite{duPCTC} and they are, for $r=1,2,3,4$ respectively $[110,111]$, $[99,103]$, $[88,97]$ and
$[77,93]$. One can see that the maximal values $\tau(r)_{max}=111,103,97$ for $r=1,2,3 =\sqrt 9$ occur only for one value of $r$, while the maximal values $\tau(4)_{max}=93$ can be a Tjurina number for curves with both $r=3$ and $r=4$. In particular, a  line arrangement with $d=12$ and $\tau(\A)\geq 94$ is  free if it has the intersection lattice of a free arrangement. As an example, the line arrangements
$$\A: f=xyz(x^9-y^9)=0 \text{ and } \A' :f'=xyz(x+y+z)(x^8-y^8)$$
are free, and have  $\tau(\A)=111$ and respectively $\tau(\A')=103$.
On the other hand, the line arrangement
$$\A: f=xyz(x^3-y^3)(y^3-z^3)(x^3-z^3)=0,$$
is free, but has $\tau(\A)=93$. Hence our result covers only a part of the free line arrangements.

\end{ex}

\begin{rk}
\label{expo}
The exponents of a free or nearly free curve of degree $d$, in particular $r=mdr(f)=d_1$, may take all the obvious possible values, see \cite{DStexpo} for examples involving both irreducible curves and line arrangements.
\end{rk}

\section{Proof of the main results} 
\subsection{Proof of Theorem \ref{thm1} and of  Corollary \ref{cor1}} 

Though this proof follows the same idea as the proof of Corollary \ref{corCTC} given above, it is longer and more involved, so we divide it into several simpler steps.

\begin{lem}
\label{lem1} If $C:f=0$ is a reduced plane curve of degree $d$ with global Tjurina number $\tau(C)=(d-1)(d-r-1)+r^2-1$, then 
$$ar(f)_{d-r-1}= {d-2r+1 \choose 2} \text{ and } n(f)_{2d-r-2}=0.$$

\end{lem}

\proof
In the case $2r<d$, we get exactly as in the subsection \ref{proofCTC} 
$$ar(f)_{d-r-1}-n(f)_{2d-r-2}={d-2r+1 \choose 2}.$$
Suppose $n(f)_{2d-r-2}>0$. Then one can reason exactly as in the proof above and obtain that $C$ is a free curve. But this is impossible, since a free curve has a different global Tjurina number by Corollary \ref{corCTC}. Hence $n(f)_{2d-r-2}=0$, which completes the proof of this Lemma when $2r<d$. If $d=2r$, then $ar(f)_{d-r-1}=ar(f)_{r-1}=0$ since $r$ is the minimal degree of a syzygy.
This clearly implies $n(f)_{2d-r-2}=0$ in this case as well.
\endproof
It was shown in \cite{DPop} that the graded $S$-module  $N(f)$ satisfies a Lefschetz type property with respect to multiplication by generic linear forms. This implies in particular the inequalities
$$0 \leq n(f)_0 \leq n(f)_1 \leq ...\leq n(f)_{[T/2]} \geq n(f)_{[T/2]+1} \geq ...\geq n(f)_T \geq 0,$$
where $T=3d-6$. Note that $T/2 <2d-r-2$ since $2r \leq d$, hence Lemma \ref{lem1} implies
\begin{equation}
\label{vanishing}
 n(f)_s=0 \text{ for any integer } s \geq 2d-r-2.
\end{equation}
It follows that the formula  \eqref{chi} for  $\chi(T\langle C\rangle(k-1))$  yields the dimension $ar(f)_{k}$ for any $k \geq d-r-1,$  namely we have 
\begin{equation}
\label{ak}
 ar(f)_k=3{k+2 \choose 2} -{d+k+1 \choose 2} +\tau(C).
\end{equation}

In particular, we get after some computation
$$ar(f)_{d-r}={d-2r-2 \choose 2}+2=\dim S_{d-2r}\rho_1+2,$$
where $\rho_1 \in AR(f)$ is the syzygy of minimal degree $r$. It follows that there are two more syzygies, say $\rho_2$ and $\rho_3$ in $AR(f)$, both of degree $d-r$, such that
$$AR(f)_{d-r}=\dim S_{d-2r}\rho_1+\C\rho_2+\C\rho_3.$$
Now we recall some basic results from \cite{duPCTC}. For two elements $ \rho=(a,b,c) \in S^3$ and $ \rho'=(a',b',c') \in S^3$, thought of as vector fields on $\C^3$, we define their exterior product in the usual way, namely
$$\rho \times \rho'=(bc'-b'c,ca'-c'a,ab'-a'b) \in S^3.$$
The following result is stated in \cite{duPCTC}, and we reprove it here in a different way for the reader convenience.
\begin{lem}
\label{lem2}
(i) If $ \rho=(a,b,c) \in S^3$ is such that the common zero set of $a,b, c$ in $\C^3$ has dimension at most 1, then $\rho \times \rho'=0$ implies that there is a polynomial $h \in S$ such that $\rho'=h\cdot \rho.$

\noindent (ii) If $ \rho$ and $ \rho'$ are two elements of $AR(f)$, then there is a polynomial $h \in S$ such that $\rho \times \rho'=h\cdot (f_x,f_y,f_z)$. 

\end{lem}

This polynomial $h$ is denoted by $\rho *\rho'$ in the sequel.

\proof To an element $ \rho=(a,b,c) \in S^3$ we can associate the differential 1-form
$$\omega(\rho)=adx+bdy+cdz.$$
Then  $\rho \times \rho'=0$ is clearly equivalent to $\omega(\rho) \wedge \omega(\rho')=0$.
The first claim (i) is a consequence of the relation between the grade of the ideal $I=(a,b,c)$ (i.e. the maximal length of a regular sequence contained in $I$) and the vanishing of the cohomology of the Koszul complex $K^*(a,b,c)$, see for instance Thm. A.2.48 in \cite{Eis}. Indeed, by our assumption on $\rho$, one has $grade(I)=2$ and the Koszul complex
$K^*(a,b,c)$ is just the complex
$$ 0 \to \Omega^0 \to \Omega^1 \to \Omega^2 \to \Omega^3 \to 0,$$
where $ \Omega^k$ denotes the $S$-module of global algebraic differential $k$-forms on $\C^3$ and the morphisms are given by the wedge product by $ \omega(\rho)$. 

To prove (ii), it is enough to check by direct computation that  $ \rho \in AR(f)$ and $ \rho' \in AR(f)$ imply that $(\rho \times \rho') \times (f_x,f_y,f_z)=0$ and the apply (i). Indeed, the common zero set of $f_x,f_y,f_z$ in $\C^3$ has dimension at most 1, since the curve $C:f=0$ is reduced.

\endproof
Following \cite{duPCTC}, we consider the sequence of graded $S$-modules
\begin{equation}
\label{es}
0 \to S(-r) \xrightarrow{u}  AR(f) \xrightarrow{v}  S(r-d+1)
\end{equation}
where the first morphism is $u:h \mapsto h\cdot \rho_1$, and the second morphism is $v:\rho \mapsto \rho *\rho_1$. By Lemma \ref{lem2}, the second morphism $v$ is well defined and the sequence is exact. Let $\ell_k=v(\rho_k)$ for $k=2,3$ and note that $\ell_2$ and $\ell_3$ are two linearly independent linear forms in $S_1$.

Assume we have a second order syzygy
\begin{equation}
\label{newsyz}
B_1\rho_1+B_2\rho_2+B_3\rho_3=0,
\end{equation}
where $B_1 \in S_k$, $B_2 \in S_{k+2r-d}$, $B_3 \in S_{k+2r-d}$ for some $k \geq d-2r+1$.
Applying the above morphism $v$ to this syzygy, we get $B_2\ell_2+B_3\ell_3=0$. It follows that
there is a polynomial $B' \in S_{k+2r-d-1}$ such that $B_2=B'\ell_3$ and $B_3=-B'\ell_2$.

Some new computation using \eqref{chi} shows that
$$ar(f)_{d-r+1}=\dim S_{d-2r+1}\rho_1+5.$$
This implies that there is at least a second order syzygy of the form 
\begin{equation}
\label{newsyz2}
C_1\rho_1+C_2\rho_2+C_3\rho_3=0,
\end{equation}
where $C_1 \in S_{d-2r+1}$, $C_2 \in S_{1}$, $C_3 \in S_{1}$.
Applying the above considerations to this syzygy, we see that there is a nonzero constant $C'$ such that  $C_2=C'\ell_3$ and $C_3=-C'\ell_2$. It follows that, up to a multiplicative constant $C'$, there is a unique second order syzygy in this degree. We call it $R$ and normalize it by choosing $C'=1$.
In higher degrees, it follows from the above that any second order sygyzy involving 
$\rho_1, \rho_2$ and $\rho_3$ is a multiple $B'\cdot R$ of the syzygy $R$. If we denote by $AR(f)'$ the graded $S$-submodule of $AR(f)$ spanned by $\rho_1, \rho_2$ and $\rho_3$ , this says exactly that we have an exact sequence
\begin{equation}
\label{es2}
0 \to S(-d+r-1) \to S(-r) \oplus S(-d+r)^2 \to AR(f)' \to 0.
\end{equation}
To complete the proof, we have to show that $ar(f)'_k=\dim AR(f)'_k$ coincides to $ar(f)_k$ for all $k \geq d-r+2$. The equality $ar(f)_k=ar(f)'_k$ for $k <d-r+2$ is obvious by the construction of $AR(f)'$, and so this would imply $AR(f)=AR(f)'$, which is exactly the property that $C$ is nearly free with exponents $(r,d-r)$.
The exact sequence \eqref{es2} implies that
$$ar(f)'_k=2{k-d+r+2 \choose 2}+ {k-r+2 \choose 2}-{k-d+r+1 \choose 2}.$$
A direct computation using this formula and the formula \eqref{ak} yields the claimed equality
$ar(f)_k=ar(f)'_k$ for $k \geq d-r+2$, thus ending the proof of Theorem \ref{thm1}.

\endproof

Finally we consider Corollary \ref{cor1}. The first claim is just a part of Proposition 1.1 in \cite{duPCTC2}, to which the reader is referred for a proof.
To prove the second claim, we use Proposition 1.3 part (2)  in \cite{duPCTC2}, were it is shown that for a reduced plane curve $C:f=0$, the condition $mdr(f)=1$ implies that either
$\tau(C)=d^2-3d+3$ (which corresponds to the case $C$ free by Corollary \ref{corCTC})
or $\tau(C)=d^2-3d+2$ (which corresponds to the case $C$ nearly free by Theorem \ref{thm1}).
When $C$ is an irreducible free curve, it is shown in Thm. 2.5 (iv) in \cite{DStFD} that the smallest exponent $d_1$ satisfies $d_1>1$. This completes the proof of Corollary \ref{cor1}.

\subsection{Proof of Theorem \ref{thmCONJ} and Corollary \ref{cor21}} 

If the curve $C:f=0$ is free, then it was shown in \cite{DStFD} that $ct(f)+st(f)=T$. So here we have to prove the stronger converse claim in $(ii)$. With the notation $r=mdr(f)$, one has $ct(f) \geq d-2+r$ as noted in 
\eqref{REL}, and hence we conclude that $st(f)\leq 2d-r-3,$ and hence
\begin{equation}
\label{eq51}
m(f)_k=\tau(C),
\end{equation}
for any $k \geq 2d-r-3.$
Proposition 2 in \cite{DBull} implies that
\begin{equation}
\label{eq52}
\dim S_k/I_{f,k}=\tau(C),
\end{equation}
for any $k \geq T-ct(f).$ In particular this holds for $k \geq 2d-r-4.$ It follows by  \eqref{eq51} and \eqref{eq52} that one has
\begin{equation}
\label{eq52.5}
n(f)_k=0 ,
\end{equation}
for any $k \geq 2d-r-3.$

It follows from Theorem 3.3 in \cite{duPCTC} that for any reduced curve $C:f=0$ of degree $d$ one has
\begin{equation}
\label{eq53}
\tau(C)=\tau(r)_{max} -e(f),
\end{equation}
where $e(f) $ is an integer such that $e(f) \geq 0$ for $r\leq (d-1)/2$ and 
$$e(f)\geq {2r+2-d \choose 2},$$
for $r>(d-1)/2$.  Then, a completely similar computation to that done in subsection \ref{proofCTC}  above yields the formula
\begin{equation}
\label{eq54}
ar(f)_{d-r-1}-n(f)_{2d-r-2}={d-2r+1 \choose 2}+1 -e(f),
\end{equation}
where the binomial coefficient, given by the usual algebraic formula, can be negative in principle, i.e. when $r$ is large. However, by \eqref{eq52.5}, we have that
the left hand side is non negative since it coincides with $ar(f)_{d-r-1}$. So the right hand side should also be non negative, and the above estimates imply that this may happen only when $r \leq d/2$. Suppose first that $r\leq (d-1)/2$. Then one has
$$\dim S_{d-2r-1}={d-2r+1 \choose 2}$$
and the exact sequence \eqref{es} implies that there are two possibilities, namely

\medskip

(a) $ar(f)_{d-r-1}=\dim S_{d-2r-1}$, or

\medskip

(b) $ar(f)_{d-r-1}=\dim S_{d-2r-1}+1.$

\medskip

The case (b) occurs if and only if $C:f=0$ is a free curve, as we have seen above. If we are in case (a), then formula \eqref{eq54} implies that $e(f)=1$. Then Theorem \ref{thm1} implies that $C:f=0$ is a nearly free curve, and for nearly free curves it was shown in \cite{DStNF} that one has
$$ct(f)+st(f)=T+2,$$
which is a contradiction with our assumption. Hence only the case (b) can occur and the proof is complete for $r\leq (d-1)/2$. 

Assume now that $r=d/2$. Hence now $d=2r$ and the formula \eqref{eq54} yields again $e(f)=1$, which we have seen to be impossible.

\medskip

We give now the  proof of Corollary \ref{cor21}. The first claim $(i)$ follows directly from Theorem \ref{thmCONJ}. As for the claim $(ii)$, the fact that a nearly free curve $C:f=0$ satisfies
$ct(f)+st(f)=T+2$ is proved in \cite{DStNF}.

Consider now the converse implication in $(ii)$, i.e. suppose that we have a curve $C:f=0$ satisfying the condition $ct(f)+st(f)=T+2$. Then as in the prove above, we see that the vanishings
\eqref{eq52.5} hold for $k \geq 2d-r-2.$ The rest of the proof of Theorem \ref{thmCONJ} applies with the only modification that now the case (b) of free curves is impossible and the remaining cases have $e(f)=1$ and hence lead to nearly free curves by Theorem \ref{thm1}.

\section{A simple characterization of nearly free curves} 
Let $C:f=0$ be a reduced plane curve of degree $d$, $r=mdr(f)$ the minimal degree of a Jacobian syzygy in $AR(f)$ and choose $\rho_1 \in AR(f)$ a homogeneous syzygy realizing this minimal degree.
Denote $S \cdot \rho_1$ the graded $S$-submodule in $AR(f)$ spanned by $\rho_1$ and consider the quotient graded module $\overline{AR(f)}=AR(f)/(S \cdot \rho_1).$ Let
$$\delta (f)_k=\dim \overline{AR(f)}_k,$$
for any integer $k$.
The exact sequence \eqref{es} gives an injection $\overline{AR(f)} \to S(r-d+1)$, and in particular we get the vanishings
\begin{equation}
\label{van2}
\delta (f)_k=0
\end{equation}
for any $k <d-r-1.$
The main result of this section is the following characterization of free and nearly free curves.
The claim (i) for free curves is just a reformulation of Lemma 1.1 in \cite{ST}, so it is already known and extensively used. The claim for nearly free curves is new and we hope useful.
\begin{thm}
\label{thm2}
With the above notation, we have the following.

\noindent (i) The curve $C$ is free if and only if $\delta (f)_{d-r-1}\geq 1.$ In such a case
$2r<d$ and  $\delta (f)_{d-r-1}=1.$

\noindent (ii) The curve $C$ is nearly free if and only if
\begin{equation}
\label{cnf}
\delta (f)_{d-r-1}=0 \text{ and } \delta (f)_{d-r} \geq 2.
\end{equation}
In such a case $2r \leq d$ and  $\delta (f)_{d-r-1}=2.$
\end{thm}

\proof
To prove (i), one can refer to the proof of Lemma 1.1 in \cite{ST}, or use the fact that in this case the second morphism $v$ in the exact sequence \eqref{es} is surjective.
So we give the details of the proof only for the case (ii). If $C$ is nearly free, it follows from the definition that the conditions \eqref{cnf} are satisfied. Conversely assume that these conditions are satisfied, and let $\rho_2, \rho_3,...,\rho_m$ be syzygies of degree $d-r$ giving a basis of 
$ \overline{AR(f)}_{d-r}$ for some $m \geq 3$. The injective morphism $ \overline{AR(f)}_{d-r} \to S_1$ induced by $v$ shows that $m\leq 4$. Let $\ell_j=v(\rho_j)$ for $j=2,3,...,m$.
Then we get $m-1$ linearly independant linear forms in $S_1$. Consider the image $V$ of the morphism $v:AR(f) \to S(r-d+1)$. Then $V$ is an ideal in $S$, and one has
$$\dim V_k =\dim S_k - \epsilon,$$
for $k>>0$, where $\epsilon \in \{0,1\}$. Indeed, $\epsilon \leq 1$ since $m\geq 3$, i.e. the ideal $V$ contains an ideal spanned by two linearly independent linear forms. Moreover $\epsilon=0$ if and only if either $m=4$, or $m=3$, but the syzygies $\rho_1,\rho_2, \rho_3$ do not generate the module $AR(f)$. Hence to prove the claim (ii) we have to show that $\epsilon=1$.
For $k>>0$, the exact sequence
$$0 \to S_{k-r} \to AR(f)_k \to V_{k+r-d+1} \to 0$$
induced by the sequence \eqref{es} yields
$$ar(f)_k={k-r+2 \choose 2}+{k+r-d+3 \choose 2}-\epsilon.$$
On the other hand, we have an obvious exact sequence
$$0 \to AR(f)_{k-d+1} \to S^3_{k-d+1} \to S_k \to M(f)_k \to 0.$$
This implies
$$m(f)_k={k+2 \choose 2}-3{k-d+3 \choose 2}+ar(f)_{k-d+1}.$$
It follows by a direct computation that the Hilbert polynomial $H(M(f))$ of $M(f)$, which is just the constant term of the $m(f)_k$ regarded as a polynomial in $k$, is given by the formula
$$H(M(f))=(d-1)(d-r-1)+r^2 -\epsilon.$$
On the other hand, it is known that $H(M(f))=\tau(C)$, see \cite{CD}.
If $\epsilon=0$, then Corollary \ref{corCTC} implies that $C$ is a free curve, which is a contradiction. Hence $\epsilon=1$, and Theorem \ref{thm1} implies that $C$ is nearly free.

\endproof

\begin{rk}
\label{diff}
The difference between Definition \ref{def2} and Theorem \ref{thm2} (ii) is that in the definition we require the three syzygies $\rho_1,\rho_2, \rho_3$ to generate $AR(f)$, a condition hard to check in practice, while in the theorem we require only the existence of (at least) two syzygies in $AR(f)_{d-r}$, independent of the syzygy $\rho_1$. This latter condition is much simpler to verify.
\end{rk}

\end{document}